\documentclass[preprint,12pt]{elsarticle}

 \usepackage{graphics}
\usepackage{amssymb}
\journal{Linear Algebra Appl.}

\begin{document}

\begin{frontmatter}

%% Title, authors and addresses

%% use the tnoteref command within \title for footnotes;
%% use the tnotetext command for the associated footnote;
%% use the fnref command within \author or \address for footnotes;
%% use the fntext command for the associated footnote;
%% use the corref command within \author for corresponding author footnotes;
%% use the cortext command for the associated footnote;
%% use the ead command for the email address,
%% and the form \ead[url] for the home page:
%%
%% \title{Title\tnoteref{label1}}
%% \tnotetext[label1]{}
%% \author{Name\corref{cor1}\fnref{label2}}
%% \ead{email address}
%% \ead[url]{home page}
%% \fntext[label2]{}
%% \cortext[cor1]{}
%% \address{Address\fnref{label3}}
%% \fntext[label3]{}

\title{The expansion of real forms on the simplex and applications}
\tnotetext[t1]{This work was supported by the National Basic
Research Program of China (973 Program: 2011CB302402), China
National Natural Science Foundation (11001228), and the Fundamental
Research Funds for the Central Universities, Southwest University
for Nationalities (12NZYTH04). }

%% use optional labels to link authors explicitly to addresses:
%% \author[label1,label2]{<author name>}
%% \address[label1]{<address>}
%% \address[label2]{<address>}

\author[yy]{Yong Yao}
\ead{yaoyong@casit.ac.cn}

\author[jx]{Jia Xu\corref{cor1}}
\ead{j.jia.xu@gmail.com}

\author[yy]{Jingzhong Zhang}
\ead{zjz101@yahoo.com.cn}

\cortext[cor1]{Corresponding author.}

\address[yy]{Chengdu Institute of Computer Applications, Chinese
Academy of Sciences, Chengdu, Sichuan 610041, China}

\address[jx]{College of Computer Science and Technology,
Southwest University for Nationalities, Chengdu, Sichuan 610041,
China}

\begin{abstract}
If $n$ points $B_1,\ldots,B_n$ in the standard simplex $\Delta_n$
are affinely independent, then they can span an $(n-1)-$simplex
denoted by $\Lambda={\rm Con}(B_1,\ldots,B_n)$. Here $\Lambda$
corresponds to an $n\times n$ matrix $[\Lambda]$ whose columns are
$B_1,\ldots,B_n$. In this paper, we firstly proved that if $\Lambda$
of diameter sufficiently small contains a point $P$, and $f(P)>0\
(<0)$ for a form $f\in{\Bbb R}[X]$, then the coefficients of
$f([\Lambda] X)$ are all positive (negative). Next, as an
application of this result, a necessary and sufficient condition for
determining the real zeros on $\Delta_n$ of a system of homogeneous
algebraic equations with integral coefficients  is established.
%\ (\ ={\Bbb R}[x_1,\ldots,x_n]\ )
\end{abstract}

\begin{keyword}
%% keywords here, in the form: keyword \sep keyword
simplex \sep system of algebraic equations \sep real zeros
%% MSC codes here, in the form: \MSC code \sep code
%% or \MSC[2008] code \sep code (2000 is the default)
\MSC 13J25 \sep 13J30 \sep 16Y60 \sep 68T15 \sep 26D05
\end{keyword}

\end{frontmatter}

%%
%% Start line numbering here if you want
%%
% \linenumbers

%% main text
\section{Introduction}
The standard simplex is denoted by $\Delta_n=\{(x_1,\ldots,x_n)^T|\
x_i\geq 0,\ \sum_i x_i=1\}$. If $n$ affinely independent points
$B_1,\ldots,B_n\in \Delta_n$, then they span an $(n-1)-$simplex
denoted by $\Lambda={\rm Con}(B_1,\ldots,B_n)$, that is,
$$\Lambda ={\rm Con}(B_1,\ldots,B_n)=\{ \lambda_1B_1+\cdots+\lambda_nB_n |
\ \sum^n_{i=1}\lambda_i=1,\ \lambda_1,\ldots,\lambda_n\geq 0\}. $$

$\Lambda$ corresponds to an $n\times n$ matrix $[\Lambda]$ whose
columns are $B_1,\ldots,B_n$. Conversely, consider the non-negative
matrix $M=[b_{ij}]$, where $B_j=(b_{1j}, \ldots,b_{nj})^T$
$(j=1,\ldots,n)$. If the determinant $|M|\neq 0$, and the sum of the
elements in each column is 1 (i.e., $\sum^n_{i=1}b_{ij}=1$), then
$M$ corresponds to an $(n-1)-$simplex ${\rm Con}(B_1,\ldots,B_n)$ in
$\Delta_n$, which is briefly denoted by ${\rm Con}(M)$.

Given an $(n-1)-$simplex $\Lambda\subseteq \Delta_n$ and a real form
(homogeneous polynomial) $\ f\in {\Bbb R}[x_1,\ldots,x_n]$, we call
$f([\Lambda]X)$ the expansion of $f$ on the simplex $\Lambda$. This
simple expansion has a surprising consequence. Here we give an
example.

\indent{\textbf{Example 1}} ( see \cite{Ly3} ) Let $(x,y,z)^T\in
\Delta_3$. Prove that
\begin{eqnarray*}
f(x,y,z)=3^7(y^4z^4(y+z)^4(2x+y+z)^8+x^4z^4(x+z)^4(x+2y+z)^8\\
+x^4y^4(x+y)^4(x+y+2z)^8)-2^8(x+y+z)^8(x+y)^4(x+z)^4(y+z)^4\geq 0.
\end{eqnarray*}

To prove this, first consider the following six matrices:
\begin{eqnarray*}
\begin{array}{ccc}
M_{1}= \left[
\begin{array}{ccc}
1 & \frac{1}{2}  & \frac{1}{3} \\
0 & \frac{1}{2} & \frac{1}{3} \\
0 &  0 & \frac{1}{3}
\end{array}
\right],\ M_{2}= \left[
\begin{array}{ccc}
1 & \frac{1}{2}  & \frac{1}{3} \\
0 &  0 & \frac{1}{3}\\
0 & \frac{1}{2} & \frac{1}{3}
\end{array}
\right],\ M_{3}= \left[
\begin{array}{ccc}
0 & \frac{1}{2} & \frac{1}{3} \\
0 &  0 & \frac{1}{3}\\
1 & \frac{1}{2}  & \frac{1}{3}
\end{array}
\right],\\
\quad\\
M_{4}= \left[
\begin{array}{ccc}
0 & \frac{1}{2} & \frac{1}{3} \\
1 & \frac{1}{2}  & \frac{1}{3} \\
0 &  0 & \frac{1}{3}
\end{array}
\right],\ M_{5}= \left[
\begin{array}{ccc}
0 &  0 & \frac{1}{3}\\
1 & \frac{1}{2}  & \frac{1}{3} \\
0 & \frac{1}{2} & \frac{1}{3}
\end{array}
\right],\ M_{6}=\left[
\begin{array}{ccc}
0 &  0 & \frac{1}{3}\\
0 & \frac{1}{2} & \frac{1}{3} \\
1 & \frac{1}{2}  & \frac{1}{3}
\end{array}
\right].
\end{array}
\end{eqnarray*}

It is easy to see that each $M_i$ corresponds to a $2-$simplex
$\Lambda_i$ in the following Figure 1 (e.g., $M_1$ corresponds to
$\Lambda_1$ and $M_2$ corresponds to $\Lambda_2$).

\begin{figure}[!htbp]
  % Requires \usepackage{graphicx}
  \centering
  \includegraphics[width=5cm]{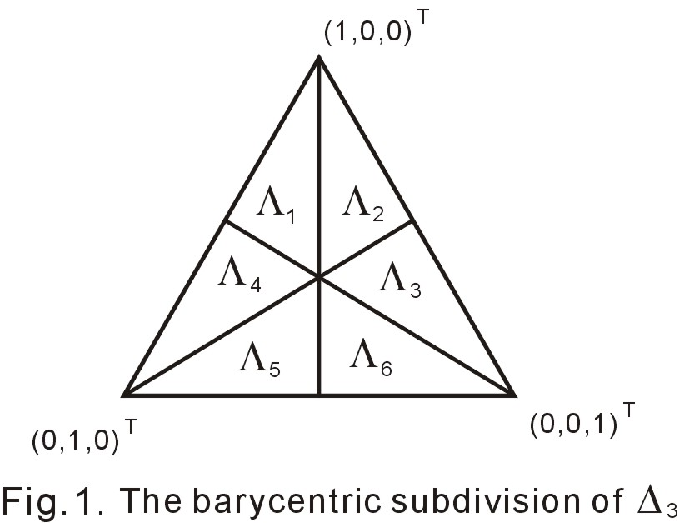}
\end{figure}

Hence we have that

\begin{eqnarray*}
\begin{array}{l}
(\forall X\in \Delta_3)\ f(X)\geq 0
 \Longleftrightarrow
 (\forall X\in \Lambda_i) \ f(X)\geq 0,\ i=1,\ldots,6, \\
  \Longleftrightarrow
 (\forall X\in \Delta_3)\ f(M_i X)\geq 0,\ i=1,\ldots,6.
\end{array}
\end{eqnarray*}

Expand the form $f(M_iX)$ with Maple (or Mathematica), we will find
that all the coefficients of $f(M_iX)$ are real non-negative
numbers. This shows that $(\forall X\in \Delta_3)\ f(X)\geq 0$.

From the above example, we know that $f([\Lambda_i]X)\
(i=1,\ldots,6)$, the expansion of $f$ on the simplex $\Lambda_i$,
can help us to determine the non-negativity of $f$ on the standard
simplex $\Delta_3$. This kind of expansion also has other
applications. This work will focus on the expansion of real forms
and applications.

Before further illustrating the main results, we will need some
notations.

Let $A=(a_1,\ldots,a_n)^T$ and $B=(b_1,\ldots,b_n)^T$ be two points
of $\Delta_n$. Then the $\infty$-norm distance $d_{\infty}(A,B)$ is
$$d_{\infty}(A,B)=\max_{1\leq i \leq n} \{|a_i-b_i|\}.$$
Furthermore, for $(n-1)-$simplex $\Lambda={\rm
Con}(B_1,\ldots,B_n)$, the diameter $\phi(\Lambda)$ is
$$\phi(\Lambda)=\max_{A,\ B\in \Lambda }\{d_{\infty}(A,\ B)\}.$$
It is well known that (see \cite{Ma1,Jr1})
$$\phi(\Lambda)=\max_{1\leq s,\ t\leq n}\{d_{\infty}(B_s,\ B_t)\},\quad  \phi(\Delta_n)=1 .$$

${\Bbb N}$ and ${\Bbb R}$ are respectively used to refer to the set
of all natural numbers and the set of all real numbers. Given  $
{\mathbf \alpha}=(\alpha_1,\ldots,\alpha_n)\in {\Bbb N}^n$, we set
$|\alpha|=\alpha_1+\cdots+\alpha_n$ and ${\mathbf
\alpha!}=\alpha_{1}!\cdots \alpha_{n}!$, and for a vector
$X=(x_1,\ldots,x_n)^T$, we have $X^{\alpha}=x^{\alpha_1}_1\cdots
x^{\alpha_n}_n.$ A form (homogeneous polynomial) $f$ of degree $d$
in ${\Bbb R^{n}}$ is given by
\begin{equation}\label{eq1}
f(X)=\sum_{|\beta|=d}\left(\frac{d!}{\beta
!}\right)C_{\beta}X^{\beta},\ \beta\in {\Bbb N}^n, \ C_{\beta}\in
{\Bbb R}.
\end{equation}
The minimum value of $f(X)$ over $\Delta_n$, the maximum value of
the coefficients of absolute value $|C_{\beta}|$ and the combination
are denoted by $\lambda_f,\ L_f$, and $C^m_{n}$ respectively.
Namely,
\begin{equation}\label{eq2}
\lambda_f=\min_{X\in \Delta_n}\{f(X)\},\ L_f=\max_{|\beta|=d}
\{|C_{\beta}|\},\ C^m_{n}=\frac{n!}{m!(n-m)!}.
\end{equation}

The main results of this paper is the following two theorems.

\newtheorem{thm1}{Theorem}
\begin{thm1}\upshape\label{thm1}
Let $f$ be in the form of (\ref{eq1}), and let $L_f, \ C^m_{n}$ be
in the form of (\ref{eq2}). The diameter of $(n-1)-$simplex $\Lambda
(\subseteq \Delta_n)$ is
 denoted by $\phi(\Lambda)$. For a point $P\in \Delta_n$ satisfying
$f(P)>0 \
 (<0)$, if $(n-1)-$simplex $\Lambda
(\subseteq \Delta_n)$ contains the point $P$ and $\phi(\Lambda)$
satisfies
$$\quad \phi(\Lambda)< (d d!C^d_{n+d-1}L_f)^{-1}|f(P)|,$$
then the coefficients of $f([\Lambda] X)$ are all positive
(negative).
\end{thm1}

\newtheorem{thm2}[thm1]{Theorem}
\begin{thm2}\upshape\label{thm2}
Let $Z=\{X\in \Delta_n |\ f_1=\cdots =f_k=0\}$, where
 $f_1,\ldots,f_k\in {\Bbb Z}[x_1,\ldots,x_n]$ are all homogeneous polynomials with degree
 $d_1,\ldots, d_k$ respectively. Denote $d=\max\{d_1,\ldots,d_k\}$. Let
 $$F=f^2_1(\sum^n_{i=1}x_i)^{2(d-d_1)}+\cdots+f^2_k(\sum^n_{i=1}x_i)^{2(d-d_k)},$$
 where $H$ is the maximum of absolute values of coefficients of $F$, and $\widetilde{H}= {\rm max} \{H,\ 4n+2\}$.
 Then the set $Z$ is not empty if and only if there is an $(n-1)-$simplex $\Lambda \subset \Delta_n$
  such that the coefficients of the polynomial $F([\Lambda]X)$ are not all positive whenever the diameter of $\Lambda$ satisfies
 $$ \phi(\Lambda)< (2d (2d)!C^{2d}_{n+2d-1}L_F)^{-1}(2^{4-\frac{n}{2}} \widetilde{H} (2d)^n)^{-n 2^n (2d)^n}.$$ Here $L_F$ is a notation in (2).
\end{thm2}
Note: In Theorem 1, "coefficients are all positive" means that every
monomial of degree $d$ appears with a strictly positive
coefficients. For example, the coefficients of the form $x^2+xy$ are
not all positive, since the coefficient of the monomial $y^2$ is
$0$.

\section{Proof of Theorem 1}

To get started, let us look at the following lemma.
\newtheorem{lem2}{Lemma}[section]
\begin{lem2}\upshape\label{lem2}
Let $f$ be in the form of (\ref{eq1}), and let $L_f$ be in the form
of (\ref{eq2}). For ${\bf \alpha}=(\alpha_1,\ldots,\alpha_n)\in
{\Bbb N}^n$, the partial derivative of $f$ is denoted by
$D^{\alpha}f$, that is,
$$
D^{\alpha}f=\frac{\partial^{\alpha_1+\alpha_2+\cdots+\alpha_n}f }
{\partial x_1^{\alpha_1}\cdots \partial x_{n}^{\alpha_n}}.
$$
Then we have
$$|D^{\alpha}f(X)|\leq d! L_{f}\ (X\in \Delta_n, \  |\alpha|\leq d).$$
\end{lem2}

\newproof{pf2}{Proof}
\begin{pf2}
Construct the polynomial
$$F(X)=L_f(x_1+\cdots+x_n)^d-f(X).$$ Since the coefficients of
$F(X)$ are all non-negative, then the coefficients of
$D^{\alpha}F(X)\ (|\alpha|\leq d)$, too, are all non-negative, that
is,
$$D^{\alpha}F(X)=D^{\alpha}(L_f(x_1+\cdots+x_n)^d)-D^{\alpha}f(X)\geq 0,\ X\in \Delta_n.$$
Hence
$$|D^{\alpha}f(X)|\leq d! L_{f}(x_1+\cdots+x_n)^{d-|\alpha|}= d! L_{f}\
(X\in \Delta_n, \  |\alpha|\leq d) .\qquad \Box$$
\end{pf2}

Next we will give the proof of Theorem 1.

\newproof{pf}{Proof of Theorem \ref{thm1}}
\begin{pf}
The basic strategy is to use the complete Taylor's formula.

Suppose that $(n-1)-$simplex $\Lambda={\rm Con}(B_1,\ldots,B_n)={\rm
Con}([b_{ij}])\subseteq \Delta_n,\ 1\leq i,j\leq n$. Let $f(P)>0 $
for the point $P=(p_1,\ldots,\ p_n)^T\in \Lambda$, and let
$\varepsilon_{ij}=b_{ij}-p_i, \
\varepsilon_{max}=\max\{|\varepsilon_{ij}|, \ 1\leq i,j\leq n \}$.
By the definition of $\phi(\Lambda)$, we know that the $\infty$-norm
distance from the point $P$ to the vertices $B_1, \ldots, B_n$ is
less than or equal to the diameter $\phi(\Lambda)$, that is,
$$|\varepsilon_{ij}|=|b_{ij}-p_i|\leq d_{\infty}(B_i,P)\leq \phi(\Lambda).$$
In particular, $\varepsilon_{max}\leq \phi(\Lambda)$.

From $b_{ij}=p_i+\varepsilon_{ij}$, it follows that
\begin{eqnarray*}
[b_{ij}]X =\left [
\begin{array}{c}
p_1|X|+(\varepsilon_{11}x_1+\cdots+\varepsilon_{1n}x_n)\\
p_2|X|+(\varepsilon_{21}x_1+\cdots+\varepsilon_{2n}x_n)\\
\cdots \\
p_n|X|+(\varepsilon_{n1}x_1+\cdots+\varepsilon_{nn}x_n)
\end{array}
\right ]=P|X|+H,
\end{eqnarray*}
where $H=(h_1,\ldots,h_n)^T, \
h_i=(\varepsilon_{i1}x_1+\cdots+\varepsilon_{in}x_n)\ and \ |X|=x_1+\cdots+x_n .$

Consider Taylor's formula
\begin{equation}\label{eqtay}
f(P|X|+H)=f(P|X|)+\sum^d_{k=1}\sum_{|\alpha|=k}\frac{D^{\alpha}
 f(P|X|)}{\alpha!}H^{\alpha},\ ({\bf \alpha}\in {\Bbb N}^n).
\end{equation}
Then we know that (\ref{eqtay}) is a identity formula for the
polynomial $f$ of degree $d$.

Notice that $f$ is homogeneous, hence by (\ref{eqtay}) we have
\begin{eqnarray*}
&f([\Lambda] X)&
=f(P|X|+H)\\
&&=f(P)|X|^d+\sum^d_{k=1}\sum_{|\alpha|=k}\frac{D^{\alpha}
 f(P)|X|^{d-k}}{\alpha!}H^{\alpha}\\
&&=\sum^d_{k=1}\left(\frac{f(P)}{d}|X|^d+\sum_{|\alpha|=k}\frac{D^{\alpha}
 f(P)}{\alpha!}|X|^{d-k}H^{\alpha}\right)\\
&&\stackrel{\star}{=} \sum^d_{k=1}\sum_{|\alpha|=k}\left(\frac{f(P)}
{dC^k_{n+k-1}}|X|^d+\frac{D^{\alpha}
 f(P)}{\alpha!}|X|^{d-k}H^{\alpha}\right)\\
&&=\sum^d_{k=1}\sum_{|\alpha|=k}
\left(\frac{f(P)}{dC^k_{n+k-1}}-\frac{|D^{\alpha}f(P)|}{\alpha!}\varepsilon_{max}
^k\right)|X|^{d}\\
&&\quad
+\sum^d_{k=1}\sum_{|\alpha|=k}\frac{|D^{\alpha}f(P)|}{\alpha!}\left(\varepsilon_{max}
^k|X|^k+H^{\alpha}\right)|X|^{d-k}.
\end{eqnarray*}
(Note: the identity $\sum_{|\alpha|=k}1=C^k_{n+k-1}$ is used at the
above sign $\stackrel{\star}{=}$.)

It remains to consider the above expression. We will do it in the
following two steps.

First, by lemma \ref{lem2}, we have
$$|D^{\alpha}f(P)|\leq d!L_f.$$
Since $$\varepsilon_{max}\leq \phi(\Lambda), \ C^k_{n+k-1}\leq
C^d_{n+d-1}\ (k\leq d).$$ Thus, we have
\begin{eqnarray*}
\left(\frac{f(P)}{dC^k_{n+k-1}}-\frac{|D^{\alpha}f(P)|}{\alpha!}\varepsilon_{max}
^k\right) \geq \left(\frac{f(P)}{dC^d_{n+d-1}}-d! L_f \phi(\Lambda)
\right) > 0.
\end{eqnarray*}

Further, it is obvious that the coefficients of
$$\frac{|D^{\alpha}f(P)|}{\alpha!}\left(\varepsilon_{max}
^k|X|^k+H^{\alpha}\right)|X|^{d-k}$$ are all non-negative real
numbers.

As stated above, we see that the coefficients of $f([\Lambda] X)$
are all positive. $\Box$
\end{pf}

\newtheorem{rmk1}{Remark}
\begin{rmk1}
For the case of $f(P)<0$, we just need to make a transformation for
$f([\Lambda] X)$:
\begin{eqnarray*}
&f([\Lambda] X)&
=f(P|X|+H)\\
&&=\sum^d_{k=1}\sum_{|\alpha|=k}
\left(\frac{f(P)}{dC^k_{n+k-1}}+\frac{|D^{\alpha}f(P)|}{\alpha!}\varepsilon_{max}
^k\right)|X|^{d}\\
&&\quad
-\sum^d_{k=1}\sum_{|\alpha|=k}\frac{|D^{\alpha}f(P)|}{\alpha!}\left(\varepsilon_{max}
^k|X|^k-H^{\alpha}\right)|X|^{d-k}.
\end{eqnarray*}
Then follow the analogous way of the above proof, we will see that
the coefficients of $f([\Lambda] X)$ are all negative when $f(P)<0$.
\end{rmk1}

Then we present a useful corollary of the theorem \ref{thm1}, which
is the basis of the successive difference substitution method (see
\cite{Ly2,Ly3,Ly4,Yy1}).

Consider the following $n\times n$ matrix $G_n$ (see \cite{Ly4,Yy1})
\begin{eqnarray*}
G_{n}= \left (
\begin{array}{cccc}
1 & \frac{1}{2} & \cdots & \frac{1}{n} \\
& \frac{1}{2} & \ddots & \vdots \\
&   &\ddots &\frac{1}{n} \\
0 &   &  & \frac{1}{n}
\end{array}
\right ).
\end{eqnarray*}
We say that a matrix is the barycentric matrix if it can got by
permuting the rows of the matrix $G_n$ (when n=3, see Example 1).

Suppose that $S_n$ is a symmetric group (permutation group) on the
set $\{1,2,\cdots,n\}$. Let $P_{\sigma}$ be a permuting matrix
corresponding to the permutation $\sigma(\sigma\in S_n)$, and let
$I$ be an identity permutation. Then the barycentric matrix can be
written as
$$G_{\sigma}=P_{\sigma}G_n=P_{\sigma}G_I, \quad \sigma\in S_n.$$
It is easy to see that the number of the barycentric matrices of
order $n$ is $n!$.

\newtheorem{cor1}{Corollary}
\begin{cor1}\upshape
Let $\ f\in {\Bbb R}[x_1,\ldots,x_n]$ be a form of degree $d$ and let $\lambda_f$ be the minimum value of $f$ on $\Delta_n$. Then:

1.\ $\lambda_f >0$ if and only if there is a positive integer $N$ such that the coefficients of
the form $f(G_{\sigma_1}G_{\sigma_2}\cdots G_{\sigma_N}X) \ (\forall \sigma_i\in S_n,\ i=1,\ldots, N)$
are all positive.

2.\ $\lambda_f <0$ if and only if there are a positive integer $N$ and a simplex ${\rm Con}(G_{\sigma_1}G_{\sigma_2}\cdots
G_{\sigma_N})$ such that
$$f(G_{\sigma_1}G_{\sigma_2}\cdots G_{\sigma_N}(1,\ldots, 1)^T)<0.$$
\end{cor1}

\newproof{pf1}{Proof}
\begin{pf1}(Sketch) To prove this corollary, we need the following three results:

(a)  For an arbitrary natural number $m$, we have that
 $$
\Delta_n=\bigcup \limits_{\sigma_1 \in
S_n} \cdots \bigcup \limits_{\sigma_m \in
S_n} {\rm Con}(G_{\sigma_1}\cdots G_{\sigma_m}).
$$

(b)  The diameter of the simplex ${\rm Con}(G_{\sigma_1}\cdots G_{\sigma_m})$ satisfies
$$
\phi ({\rm Con}(G_{\sigma_1}\cdots G_{\sigma_m}))\leq \left(\frac{n-1}{n}\right )^m.
$$

(c)  There is a positive integer $N$ satisfying
\begin{equation}\label{eq3}
\left(\frac{n-1}{n}\right)^N<(d d!
C^d_{n+d-1}L_f)^{-1}|\lambda_f|.
\end{equation}

Here (a) and (b) are the immediate consequence of barycentric
subdivision of the simplex $\Delta_n$ (see \cite{Ma1,Jr1}). For a
given form $f$, the right hand of the inequality (\ref{eq3}) is a
constant whenever its left hand can be arbitrarily close to $0$.
Hence such $N$ that satisfies the (\ref{eq3}) is existent. Then by
applying Theorem 1, the corollary follows immediately. $\Box $
\end{pf1}

\section{Proof of Theorem 2}
The discussion in this section is restricted in the polynomial ring
${\Bbb Z}[x_1,\ldots,\allowbreak x_n]$. The following lemma is
essential in proving Theorem 2.

\newtheorem{lem}{Lemma}[section]
\begin{lem}\upshape
(see \cite{Gj1}) Let $T = \{x\in {\Bbb R}^n\ | f_1(x)= \ldots =
f_l(x)=0, f_{l+1}(x)\allowbreak \geq 0, \ldots, f_m(x) \geq 0\}$ be
defined by polynomials $f_1, \ldots, f_m \in {\Bbb
Z}[x_1,\ldots,x_n]$ with degrees bounded by an even integer $d$ and
coefficients of absolute value at most $H$, and let $C$ be a compact
connected component of $T$. Let $g \in {\Bbb Z}[x_1,\ldots,x_n]$ be
a polynomial of degree $d_0\leq d$ and coefficients of absolute
value bounded by $H_0\leq H$. Then, the minimum value that $g$ takes
over $C$ is a real algebraic number of degree at most
$$2^{n-1}d^n$$
and, if it is not zero, its absolute value is greater or equal to
$$(2^{4-\frac{n}{2}}\widetilde{H} d^n)^{-n 2^{n}d^{n}},$$
where $\widetilde{H}= {\rm max} \{H,2n+2m\}$.
\end{lem}
\newproof{pf3}{Proof of Theorem \ref{thm2}}
\begin{pf3}
We first show that if the coefficients of $F([\Lambda]X)$ are not
all positive, then the set $Z$ is not empty. We will prove this by
contradiction. It is clear that $Z=\{X\in \Delta_n | \ F=0 \}$ since
$Z=\{X\in \Delta_n |\ f_1=\cdots =f_k=0\}$. Let
$T=\Delta_n=\{(x_1,\ldots,x_n)^T|\ \sum_i x_i=1,\ x_i\geq 0 \}$, and
let $g=F$.
 By Lemma 3.1, if there is no root of $F$ on $T$, then the minimum value $\lambda_F$
 of $F$ on $T$ satisfies $$\lambda_F\geq (2^{4-\frac{n}{2}}\widetilde{H} (2d)^n)^{-n 2^{n}(2d)^{n}},$$
where $\widetilde{H}= {\rm max} \{H,\ 4n+2\}$.

Thus, for a simplex $\Lambda$, the diameter $\phi(\Lambda)$
satisfies

$$\phi(\Lambda)<(2d (2d)!C^{2d}_{n+2d-1}L_F)^{-1}\lambda_F,$$
since it satisfies
$$ \phi(\Lambda)< (2d (2d)!C^{2d}_{n+2d-1}L_F)^{-1}(2^{4-\frac{n}{2}} \widetilde{H} (2d)^n)^{-n 2^n (2d)^n}.$$

 By Theorem \ref{thm1}, we know that the coefficients of the form $F([\Lambda]X)$ are all positive,
  which contradicts the condition that the coefficients of $F([\Lambda]X)$ are not all positive.
  Hence $Z$ is not empty. The converse is trivial and the theorem is proved. $\Box$
\end{pf3}

Finally, we should remark that if you are interested in the other
applications about the expansion of real forms on the simplex,
please see \cite{Jx1,Jx2,Jx3,Ly2,Ly3,Ly4,Yy1}.

%% The Appendices part is started with the command \appendix;
%% appendix sections are then done as normal sections
%% \appendix

%% \section{}
%% \label{}

%% References
%%
%% Following citation commands can be used in the body text:
%% Usage of \cite is as follows:
%%   \cite{key}          ==>>  [#]
%%   \cite[chap. 2]{key} ==>>  [#, chap. 2]
%%   \citet{key}         ==>>  Author [#]

%% References with bibTeX database:

\bibliographystyle{model1-num-names}
%\bibliography{<your-bib-database>}

%% Authors are advised to submit their bibtex database files. They are
%% requested to list a bibtex style file in the manuscript if they do
%% not want to use model1-num-names.bst.

%% References without bibTeX database:

%% \bibitem must have the following form:
%%   \bibitem{key}...
%%

%% \bibitem{}

\end{document}